\input amstex
\documentstyle{amsppt}
\TagsOnRight
\magnification=\magstep1
\pageheight{24 truecm}
\pagewidth{16.8 truecm}
\input diagrams
%
\loadbold
\loadeusm
\def\scr#1{{\fam\eusmfam\relax#1}}

%
%
\def\openC{\Bbb C}

\def\openP{\Bbb P}

\def\openZ{\Bbb Z}

\def\Hom{\operatorname{Hom}}

\def\dim{\operatorname{dim}}
\def\ker{\operatorname{ker}}

\def\Hom{\operatorname{Hom}}
\def\Ext{\operatorname{Ext}}

\def\red{\operatorname{red}}
\def\C{\operatorname{C}}

\def\rto{\raise.5ex\hbox{$\scriptscriptstyle ---\!\!\!>$}}
%
%
%
\topmatter
\title
An elliptic conic bundle in ${\openP}^4$ arising
from a stable rank-3 vector bundle
\endtitle
\rightheadtext{H. Abo, W. Decker, N. Sasakura}
\leftheadtext{An elliptic conic bundle in ${\openP}^4$}

\author
Hirotachi Abo, Wolfram Decker, Nobuo Sasakura\dag
\endauthor
%
%
\address
Hirotachi Abo\hfil\break
Department of Mathematics, 
Tokyo Metropolitan University,
Minami-Ohsawa 1-1, Hachioji-shi
Tokyo, 192-03
Japan 
\endaddress
\email
abo\@math.metro-u.ac.jp
\endemail
\address
Wolfram Decker\hfil\break
Fachbereich Mathematik,
Universit\"at des Saarlandes,
D 66041 Saarbr\"ucken,
Germany
\endaddress
\email
decker\@math.uni-sb.de
\endemail
%
%
%
\endtopmatter
\document
%
\head
{0. Introduction}
\endhead
In this paper we give a cohomological proof for the existence 
of a new family of smooth surfaces in $\openP^4$. In fact, we show:
\proclaim{Theorem 0.1} 
(i)\quad Let $X\subset\openP^4$ be an elliptic conic bundle with
degree $d=8$ and sectional genus $\pi =5$. Then the ideal sheaf of 
$X$ arises as a cokernel
$$
0\rightarrow 4\scr O(-1)\rightarrow \scr G \rightarrow \scr J_X(3)
\rightarrow 0, 
$$
where $\scr G$ is a rank-5 vector bundle on $\openP^4$ with Chern-classes
$$
c_1=-1, c_2=2, c_3=-2 \quad {\text {and}} \quad c_4=-3.
$$
$\scr G$ is isomorphic to the cohomology bundle of a monad
$$
(M) \qquad 0\rightarrow \Omega ^3(3) \overset\alpha\to\rightarrow
\Omega ^2(2) \oplus \Omega ^1(1) \overset\beta\to\rightarrow 
\scr O\rightarrow 0,
$$
with
$$
\alpha=
\pmatrix
e_4\\
e_0\wedge e_2 + e_1\wedge e_3
\endpmatrix
$$
and
$$
\beta=
\pmatrix
    &e_0\wedge e_2 + e_1\wedge e_3&-e_4\\
\endpmatrix ,
$$
\smallskip\noindent
where $e_0,\dots,e_4$ is a basis of the underlying vector space $V$
of $\openP^4$. In particular, $\scr G$ is uniquely determined up to
isomorphisms and coordinate transformations.\par\noindent
(ii)\quad Conversely, if $\scr G$ is the cohomology bunndle of the
monad $(M)$ as in (i), then $\scr G(1)$ is globally generated. Therefore 
the dependancy locus of four general sections of $\scr G(1)$ is a
smooth surface $X\subset\openP^4$. In fact, $X$ is an
elliptic conic bundle with invariants as above.\qed
\endproclaim
\noindent
For a geometric construction of our surfaces we refer to \cite{Ra}.\medskip
The new surfaces are missing in a series of classification papers. First 
of all, they are falsely ruled out in the classification of smooth degree 8 
surfaces in $\openP^4$ by Okonek \cite {Ok2} and, independently, Ionescu 
\cite {Io}.
As a consequence, they are e.g. also missing in the first version
of two papers which are concerned with the classification of conic bundles
in $\openP^4$ \cite {ES}, \cite {BR}. The correct result is:
\proclaim{Theorem 0.2} \cite{ES}, \cite {BR}
Let $X\subset \openP^4$ be a smooth surface ruled in conics. Then $X$ is
either rational or an elliptic conic bundle with $d=8$ and $\pi =5$. In 
the first case $X$ is either a Del Pezzo surface of degree 4 or a Castelnuovo
surface.\qed
\endproclaim
\noindent
This result is important in the context of adjunction theory \cite {So}, 
\cite {SV}, \cite {VdV}. Recall, that for a smooth surface $X\subset\openP^4$
the adjunction map is defined unless $X$ is a plane or a scroll. If the 
adjunction map is defined, then it has a 2-dimensional image unless $X$ is a 
Del Pezzo surface or a conic bundle. Therefore the classification of scrolls 
in $\openP^4$ \cite{La}, \cite {Au} and the above result imply:
\proclaim{Corollary 0.3}
Let $X\subset \openP^4$ be a smooth surface of degree $d\geq 9$. Then the
adjunction map is defined and has a 2-dimensional image. \qed
\endproclaim
\noindent
The new family is one of a few known families of irregular smooth surfaces 
in $\openP^4$. In fact, up to pullbacks via finite covers 
$\openP^4\to\openP^4$, our surfaces are the first such surfaces which do not 
possess a Heisenberg symmetry (compare \cite {ADHPR}). Moreover, they provide 
a counterexample to a conjecture of 
Ellingsrud and Peskine. According to this conjecture there should be no irregular 
m-ruled surface in $\openP^4$ for $m\geq 2$.\par
\medskip
We first came across the elliptic conic bundles when studying a stable 
rank-3 vector bundle $\scr E$ on $\openP^4$ with Chern classes
$$
c_1=4, c_2=8 \quad {\text {and}} \quad  c_3=8.
$$
$\scr E$ has been found by the stratification theoretical method of the third 
author (compare \cite{Sa} for this method). The dependancy locus of two 
sections of $\scr E$ is a smooth surface of the desired type. In fact, 
$\scr E$ is a cokernel
$$
0\rightarrow 2 \scr O\rightarrow\scr G(1) \rightarrow \scr E\rightarrow 0,
$$
where $\scr G$ is the rank-5 vector bundle from Theorem 0.1.\medskip
Our paper is organized as follows:
In Section 1 we review Beilinson's theorem \cite{Bei} in the context of smooth
surfaces in $\openP^4$. In Section 2 we follow the cohomological approach 
of Okonek and find the rank-5 vector bundle $\scr G$ and thus the elliptic
conic bundles via Beilinson's generalized monads. In Section 3 we take a quick
look at the stratification theoretical method of the third author
by studying some examples. In particular, we will construct the rank-3 bundle
$\scr E$ and we will explain its relation to the rank-5 bundle $\scr G$.\medskip
The third author, Nobuo Sasakura, died in June 1997. His death has 
been a great loss to us.\medskip
\noindent
{\it {Acknowledgement}}. We thank Kristian Ranestad and Frank-Olaf Schreyer
for helpful discussions. The second author is grateful to the Tokyo Metropolitan 
University and to M. Maruyama and the Kyoto University for their hospitality.
The first author would like to thank the German research society 
for its financial support (DFG-Schwerpunkt "Algorith\-mi\-sche Zahlentheorie und
Algebra") and the Universit\"at des Saarlandes for its hospitality. 
\par
\proclaim{Notations 0.4} $\openP^4 = \openP^4(x) =\openP (V)$ will be the 
projective space of lines in a $5$-dimensional $\openC$-vector space $V$ and 
$$
R = \openC [x_0,\dots, x_4] = \underset {m\geq 0}\to{\textstyle\bigoplus} 
S^m V^*
$$
its homogeneous coordinate ring. We write $\Omega^i = \Lambda^i
\scr T_{\openP^4}^*$ for the $i^{\text {th}}$ bundle of differentials
on $\openP^4$. If $X$ is a smooth surface in $\openP^4$, then we denote by\par
\noindent -\quad $H$ its hyperplane class\par
\noindent -\quad $d=H^2$ its degree\par
\noindent -\quad $K=K_X$ its canonical divisor\par
\noindent -\quad $\pi =\frac 1 2 H\!\centerdot\! (H+K)+1$ its sectional 
genus\par
\noindent -\quad $q=p_g-p_a$ its irregularity\par 
\noindent -\quad $\chi=\chi(\scr O_X)=1-q+p_g$ the Euler characteristic
of its structure sheaf.\quad\qed
\endproclaim
\bigskip\noindent
\head
{1. Beilinson's theorem}
\endhead
In this section we review Beilinson's theorem \cite{Bei} in the context of smooth
surfaces in $\openP^4$. The general theorem tells us that the derived 
category of coherent sheaves on $\openP^n$ is generated by the twisted bundles 
of differentials $\Omega^i(i),$ $0\leq i \leq n$. More concretely, every coherent 
sheaf on $\openP^n$ is the cohomology of a certain generalized monad 
involving bundles of differentials. In the next section we apply Beilinson's 
theorem in order to study our surfaces in $\openP^4$ from a cohomological 
point of view.
\proclaim{Theorem 1.1}\cite{Bei}\quad
For any coherent sheaf $\scr S$ on $\openP^4$ there is a complex
$\scr K^{\cdot}$ with 
$$
\scr K^i \cong \bigoplus _j H^{i+j}\left(\openP^4, \scr S(-j)\right)\otimes 
\Omega^j(j),
$$
such that
$$
H^i(\scr K^{\cdot}) = \left\{
\aligned
\scr S&\quad i=0\\
0&\quad i\neq 0\quad .\quad \qed
\endaligned
\right . 
$$
\endproclaim 
\noindent
The differentials of $\scr K^{\cdot}$ are given by matrices with entries in 
the exterior algebra $\Lambda V$ over the underlying vector space $V$
of $\openP^4$:
\proclaim
{Remark 1.2} The Koszul complex on $\openP^4=\openP (V)$ is the exact sequence
{\eightpoint
$$
0\leftarrow \scr O\overset{\kappa_1}\to\leftarrow
V^*\otimes \scr O(-1)\overset{\kappa_2}\to\leftarrow 
\Lambda ^2 V^* \otimes \scr O(-2)\overset{\kappa_3}\to\leftarrow \quad \dots 
\quad \overset{\kappa_{5}}\to\leftarrow 
\Lambda ^{5} V^*\otimes  \scr O(-5))\leftarrow 0,\ 
$$}\noindent
defined by contraction with the tautological subbundle
$$
\scr O(-1)\rightarrow V\otimes \scr O.
$$
The syzygy bundles of this complex
are the bundles of differentials:
$$
(\ker \kappa_i) \cong \Omega^i,\quad 0\leq i\leq 4.
$$
Via the short exact sequences associated to the Koszul complex we may compute 
the intermediate cohomology modules of the $\Omega^i$:
$$
H^j_{\ast} \left(\openP^4, \Omega^i\right) \cong \left\{
\aligned
\openC&\quad 1\leq j=i\leq 3\\
0&\quad 1\leq j\neq i\leq 3\quad .
\endaligned
\right .
$$
In the same way we obtain for $0\leq  i, j\leq 4$:
$$
\Hom(\Omega^i(i), \Omega^j(j))\cong\Lambda ^{i-j}V,
$$
the isomorphism being defined by contraction. The composition
of homomorphisms coincides with multiplication in $\Lambda V$.
\qed
\endproclaim
\noindent
We may therefore express conditions on a coherent sheaf on $\openP (V)$ 
as conditions on certain matrices with entries in $\Lambda V$.
For example we need:
\proclaim
{Remark 1.3} Let 
$$
t\Omega^i(i)\overset A\to\rightarrow s\Omega^{i-1}(i-1)
$$
be a vector bundle homomorphism, i.e., let $A$ be a 
$s\times t$-matrix with entries in $V$. Then a necessary condition for $A$ to 
be surjective is: If $(a_1,\dots, a_t)$ is a non-trivial
linear combination of the rows of $A$, then
$$
\dim {\def\span{\operatorname{span}} \span} (a_1,\dots, a_t)\geq i+1:
$$
$A$ is surjective iff its dual map is pointwise
injective iff
$$
s\Lambda^{i-1}V\wedge x\overset A\to\rightarrow t\Lambda^{i}V\wedge x
$$
is injective for any $\langle x\rangle\in\openP (V)$. 
\qed
\endproclaim 
\noindent
In order to apply Beilinson's theorem to the ideal sheaf
(suitably twisted) of a smooth surface $X\subset\openP^4$
we need information on the dimensions $h^i\scr J_X (m)$. First 
we recall Riemann-Roch:
\proclaim{Proposition 1.4 }
Let $X\subset\openP^4$ be a smooth surface. Then
$$
\chi(\scr J_X(m)) = \chi(\scr O_{\openP^4}(m))-\binom{m+1}2 d
+m(\pi-1)-1+q-p_g\ .\qquad\qed
$$
\endproclaim
\noindent Moreover one knows:
\proclaim{Proposition 1.5} \cite{DES}
Let $X\subset \openP^4$ be a smooth, non-general type surface which is
not contained in any cubic hypersurface. Then we have the following
table for the $h^i\scr J_X (m)$:
%
$$
\vbox{\offinterlineskip
\halign{\hbox to 15pt{\hfil}
\vrule height10.5pt depth 5.5pt#
&\hbox to 50pt{\hfil$#$\hfil}
&\vrule#
&\hbox to 50pt{\hfil$#$\hfil}
&\vrule#
&\hbox to 50pt{\hfil$#$\hfil}
&\vrule#
&\hbox to 50pt{\hfil$#$\hfil}
&\vrule#
&\hbox to 50pt{\hfil$#$\hfil}
&\vrule#
&\hbox to 50pt{\hfil$#$\hfil}
&\vrule#
&\hbox to 15pt{\hfil$#$\hfil}
&#
\cr
\omit&&\omit\hbox{\hskip-3pt\hbox{$\bigg\uparrow$}\raise5pt\hbox{$i$}\hfil}\cr
\multispan{14}\hrulefill\cr
&0&&0&&0&&0&&0&&0&& \cr
\multispan{14}\hrulefill\cr
&N+1&&p_g&&0&&0&&0&&0&\cr
\multispan{14}\hrulefill\cr
&0&&q&&h^2\scr J_X(1)&&h^2\scr J_X(2)&&h^2\scr J_X(3)&&h^2\scr J_X(4)&\cr
\multispan{14}\hrulefill\cr
&0&&0&&0&&h^1\scr J_X(2)&&h^1\scr J_X(3)&&h^1\scr J_X(4)&\cr
\multispan{14}\hrulefill\cr
&0&&0&&0&&0&&0&&h^0\scr J_X(4)&\cr
\multispan{14}\hrulefill
&\vbox to0pt{\vss\vskip5.75pt\hbox{$\!-\negthickspace\negmedspace@>> m >$}\vss}\cr
}}
$$
\par\noindent
where
$$
N = \pi - q + p_g - 1\ .\qquad\qed
$$
\endproclaim\noindent
In the sequel we represent a zero in a cohomology table
by an empty box.
\bigskip\noindent
\head
{2. A monad construction}
\endhead
In this section we show, that every elliptic conic bundle in $\openP^4$ 
with $d=8$ and $\pi=5$ is given as the dependancy locus of four sections
in a unique rank-5 vector bundle $\scr G(1)$. We give a monad construction
for $\scr G$, which conversely provides a proof for the existence of 
our surfaces since $\scr G(1)$ is globally generated. 
\medskip
We first mention a classification result concerning the numerical 
invariants of those surfaces we are interested in here. This result is a
consequence of generalized Serre correspondence 
\cite{Ok1, Theorem 2.2} and adjunction theory \cite{So}, \cite{SV},
\cite{VdV}.
\proclaim
{Proposition 2.1} \cite{Ok2}\quad Let $X\subset\openP^4$ 
be a smooth surface with $d=8$ and $\pi=5$. Then $p_g=0$ 
and $q\leq 1$. Moreover:\par\noindent 
(i)\quad If $q=0$, then $X$ is the blow-up of $\openP^2$ 
in eleven points, 
$$X=\openP^2(p_0,\dots, p_{10}),$$ 
embedded by the linear system (with obvious notations)
$$H=7L-E_0-\sum_{i=1}^{10}2E_i.$$
\par\noindent  
(ii)\quad If $q=1$, then the adjunction map $\Phi_{\mid K+H \mid}$ exhibits 
$X$ as a conic bundle over an elliptic normal curve $C$ in $\openP^3$. There
are precisely eight singular fibres of $\Phi_{\mid K+H \mid}$.
These singular fibres are pairs of (-1)-lines $E_i, 
\tilde E_i$ with $E_i\cdot \tilde E_i=1.$ 
\endproclaim
\demo
{Proof of (ii)} By classification (see \cite {La} and \cite {Au}) $X$ is
not a scroll. It follows from \cite {So}, \cite {VdV} that $\scr O_X(K+H)$ 
is spanned, i.e., that the adjunction map is a well-defined morphism
$X\rightarrow\openP^N$ with $N=\pi - q + p_g - 1=3$. We have 
$H \!\centerdot\! K=0$ 
since $\pi=5$. So we obtain $K^2=-8$ and thus $(K+H)^2=0$ from the double 
point formula (cf. \cite{Ha, Appendix A, 4.1.3})
$$
d^2 - 10d - 5 H \!\centerdot\! K -2K^2 +12 \chi = 0. 
$$
By \cite{So} $\Phi_{\mid K+H \mid}$ exhibits $X$ as a conic bundle
over a smooth elliptic curve in $\openP^3$ of degree $\pi-1=4$.
The singular fibres of $\Phi_{\mid K+H \mid}$
are pairs of (-1)-lines $E_i, \tilde E_i$ with $E_i\cdot \tilde E_i=1.$
By blowing down one irreducible component for each singular fiber
we obtain a ruled surface $Y$ of irregularity $q=1$ together with a 
commutative diagram
$$
\diagram[height=2.5em]
X&&\rTo^{\Phi_{\mid K+H \mid}}&&{C\subset\openP^3}\\
&\rdTo&&\ruTo\\
&&Y&.\\
\enddiagram
$$
By comparing $K_X^2=-8$ with $K_Y^2=8(1-q)=0$ we find that $\Phi$ has 
precisely eight singular fibers.
\quad\qed
\enddemo
\noindent 
The existence of surfaces of type (i) has been verified by Alexander
\cite {Al}. Okonek \cite{Ok2} mistakenly claimed that surfaces $X$ of
type (ii) do not exist. Essentially, he claimed that the generalized 
monad of Beilinson for the twisted ideal sheaf $\scr J_X(2)$ cannot exist.
He gave no argument, however, and in fact such a monad can be easily 
constructed. It is even more convenient to apply Beilinson's theorem 
to $\scr J_X(3)$. 
\proclaim
{Proposition 2.2} 
Let $X\subset\openP^4$ be a smooth surface with $d=8$, $\pi = 5$
and $q=1$. Then:\par
\noindent (i)\quad $\scr J_{X}$  has the following cohomology table:
$$
\vbox{\offinterlineskip
\halign{\hbox to 10pt{\hfil#}
&\vrule height10.5pt depth 5.5pt#
&\hbox to 25pt{\hfil$#$\hfil}
&\vrule#
&\hbox to 20pt{\hfil$#$\hfil}
&\vrule#
&\hbox to 25pt{\hfil$#$\hfil}
&\vrule#
&\hbox to 25pt{\hfil$#$\hfil}
&\vrule#
&\hbox to 25pt{\hfil$#$\hfil}
&\vrule#
&\hbox to 25pt{\hfil$#$\hfil}
&\vrule#
&\hbox to 10pt{#\hss}
\cr
&\omit&&\omit\hbox{\hskip-3pt\hbox{$\bigg\uparrow$}\raise5pt\hbox{$i$}\hfil}\cr
\multispan{15}\hrulefill\cr
&& && && && && && &\cr
&\multispan{12}\hrulefill\cr
&&4&& && && && && &\cr
&\multispan{12}\hrulefill\cr
&& &&1&&1&& && && &\cr
&\multispan{12}\hrulefill\cr
&& && && &&1&&1&& &\cr
&\multispan{12}\hrulefill\cr
&& && && && && &&6&\cr
\multispan{14}\hrulefill
&\vbox to 0pt{\vss\vskip 4pt\hbox to 40pt{\hskip-.5pt\rightarrowfill}
\hbox to 33pt{\hfil$m$}\vss}\cr
}}
$$
In particular, $X$ is cut out by 6 quartics.
\medskip
\noindent (ii) There exists an exact sequence
$$
0\rightarrow 4\scr O(-1)\rightarrow \scr G \rightarrow \scr J_X(3)
\rightarrow 0, 
$$
where $\scr G$ is a rank-5 vector bundle on $\openP^4$ with Chern-classes
$$
c_1=-1, c_2=2, c_3=-2 \quad {\text {and}} \quad c_4=-3.
$$
\noindent 
$\scr G$ is isomorphic to the cohomology bundle of a monad
$$
(M) \qquad 0\rightarrow \Omega ^3(3) \overset\alpha\to\rightarrow
\Omega ^2(2) \oplus \Omega ^1(1) \overset\beta\to\rightarrow 
\scr O\rightarrow 0,
$$
with
$$
\alpha=
\pmatrix
e_4\\
e_0\wedge e_2 + e_1\wedge e_3
\endpmatrix
$$
and
$$
\beta=
\pmatrix
    e_0\wedge e_2 + e_1\wedge e_3&-e_4
\endpmatrix ,
$$
\smallskip\noindent
where $e_0,\dots,e_4$ is a basis of the underlying vector space $V$
of $\openP^4$. In particular, $\scr G$ is uniquely determined up to
isomorphisms and coordinate transformations.
\medskip
\noindent (iii)\quad $\scr G$ has a minimal free resolution of type
$$
\vbox{%
\halign{&\hfil\,$#$\,\hfil\cr
0&\leftarrow&\scr G&\leftarrow&&10 \scr O(-1)&&4\scr O(-2)&&\scr O(-3)\cr
&&&&&&\vbox to 10pt{\vskip-0pt\hbox{$\nwarrow$}\vss}&\oplus&\leftarrow&\oplus\cr
&&&&&&&5\scr O(-3) &&4\scr O(-4)&\vbox to 10pt{\vskip-4pt\hbox{$\nwarrow$}\vss}&
\scr O(-5)&\leftarrow&0\ .\cr
}}
$$
So $\scr J_X$ has syzygies of type
$$
\vbox{%
\halign{&\hfil\,$#$\,\hfil\cr
0&\leftarrow&\scr J_X&\leftarrow&&6 \scr O(-4)&&4\scr O(-5)&&\scr O(-6)\cr
&&&&&&\vbox to 10pt{\vskip-0pt\hbox{$\nwarrow$}\vss}&\oplus&\leftarrow&\oplus\cr
&&&&&&&5\scr O(-6) &&4\scr O(-7)&\vbox to 10pt{\vskip-4pt\hbox{$\nwarrow$}\vss}&
\scr O(-8)&\leftarrow&0\ .\cr
}}
$$
\endproclaim
\demo
{Proof}
(i)\quad $X$ is a ruled surface by Proposition 2.1, (ii). In
particular, $X$ is not of general type. Moreover, $X$ is not
contained in a cubic hypersurface by general classification 
results (see \cite {Ro}, \cite {Au} and \cite {Ko}). Therefore 
$\scr J_X$ has a cohomology table of type
$$
\vbox{\offinterlineskip
\halign{\hbox to 10pt{\hfil#}
&\vrule height10.5pt depth 5.5pt#
&\hbox to 25pt{\hfil$#$\hfil}
&\vrule#
&\hbox to 20pt{\hfil$#$\hfil}
&\vrule#
&\hbox to 25pt{\hfil$#$\hfil}
&\vrule#
&\hbox to 25pt{\hfil$#$\hfil}
&\vrule#
&\hbox to 25pt{\hfil$#$\hfil}
&\vrule#
&\hbox to 25pt{\hfil$#$\hfil}
&\vrule#
&\hbox to 10pt{#\hss}
\cr
&\omit&&\omit\hbox{\hskip-3pt\hbox{$\bigg\uparrow$}\raise5pt\hbox{$i$}\hfil}\cr
\multispan{15}\hrulefill\cr
&& && && && && && &\cr
&\multispan{12}\hrulefill\cr
&&4&& && && && && &\cr
&\multispan{12}\hrulefill\cr
&& &&1&&1&&a&&*&&*&\cr
&\multispan{12}\hrulefill\cr
&& && && &&*&&*&&*&\cr
&\multispan{12}\hrulefill\cr
&& && && && && &&*&\cr
\multispan{14}\hrulefill
&\vbox to 0pt{\vss\vskip 4pt\hbox to 40pt{\hskip-.5pt\rightarrowfill}
\hbox to 33pt{\hfil$m$}\vss}\cr
}}
$$
by Proposition 1.5. Beilinson's theorem applied to $\scr J_X(2)$ implies
the existence of a surjective map $\Omega^1(1)\rightarrow a\scr O$.
By Remark 1.3 this is only possible if $a = h^2 \scr J_X(2)=0$. Then
$h^2 \scr J_X(m)=0$ for every $m\geq 2$. In particular, 
$h^1 \scr J_X(2)=h^1 \scr J_X(3)=1$ by Riemann-Roch. Applying the same 
argument again gives $h^1 \scr J_X(m)=0$ for every $m\geq 4$ and 
$h^0 \scr J_X(4)=6$.
\par\noindent
(ii)\quad The adjoint linear system
$$
H^0(\scr O_X(K+H))\cong \Ext^1(\scr J_X(3), \scr O(-1))=:W
$$
has dimension 4 and generates $\scr O_X(K+H)\cong\omega_X(1)$
(compare the proof of Proposition 2.1). Therefore the identy in
$$
W^*\otimes W\cong \Ext^1(\scr J_X(3), W^*\otimes \scr O(-1))
$$
defines an extension as in the assertion with a locally free sheaf $\scr G$. 
By construction and (i), $\scr G(-3)$ has the following cohomology table:
$$
\vbox{\offinterlineskip
\halign{\hbox to 10pt{\hfil#}
&\vrule height10.5pt depth 5.5pt#
&\hbox to 25pt{\hfil$#$\hfil}
&\vrule#
&\hbox to 20pt{\hfil$#$\hfil}
&\vrule#
&\hbox to 25pt{\hfil$#$\hfil}
&\vrule#
&\hbox to 25pt{\hfil$#$\hfil}
&\vrule#
&\hbox to 25pt{\hfil$#$\hfil}
&\vrule#
&\hbox to 10pt{#\hss}
\cr
&\omit&&\omit\hbox{\hskip-3pt\hbox{$\bigg\uparrow$}\raise5pt\hbox{$i$}\hfil}\cr
\multispan{13}\hrulefill\cr
&& && && && && &\cr
&\multispan{10}\hrulefill\cr
&& && && && && &\cr
&\multispan{10}\hrulefill\cr
&& &&1&&1&& && &\cr
&\multispan{10}\hrulefill\cr
&& && && &&1&&1&\cr
&\multispan{10}\hrulefill\cr
&& && && && && &\cr
\multispan{12}\hrulefill
&\vbox to 0pt{\vss\vskip 4pt\hbox to 40pt{\hskip-.5pt\rightarrowfill}
\hbox to 33pt{\hfil$m$}\vss}\cr
}}
$$
Beilinson's theorem implies that $\scr G$ is the cohomology bundle
of a monad with bundles of differentials as in the assertion. 
Expressing the monad conditions in terms of the exterior algebra
$\Lambda V$ shows, that, up to isomorphisms, also the arrows of the 
monad are of the claimed type (compare Remark 1.2 and Remark 1.3).
\par\noindent
(iii)\quad 
The kernel $\scr K=\ker \beta$ is a rank-9 bundle which fits
into a commutative diagram with exact rows and columns as follows:
{\eightpoint
$$
\diagram[height=2.5em]
&&{0}&&{0}\\
&&\dTo&&\dTo\\
{0}&\rTo&{\scr K}&\rTo&{\Omega^2(2)\oplus\Omega^1(1)}&\rTo^\beta
&{\scr O}&\rTo&{0}\\
&&\dTo&&\dTo&&||\\
{0}&\rTo&{U\otimes\scr O}&\rTo&{(\Lambda^2 V^*\oplus \Lambda^1 V^*)
\otimes\scr O}&\rTo^\beta&{\scr O}&\rTo&{0}\\
&&\dTo^\phi&&\dTo\\
&&{(\Lambda^1 V^*\oplus \Lambda^0 V^*)\otimes\scr O(1)}&=&
{(\Lambda^1 V^*\oplus \Lambda^0 V^*)\otimes\scr O(1)}\\
&&\dTo&&\dTo\\
&&{\scr O(2)}&=&{\scr O(2)}\\
&&\dTo&&\dTo\\
&&{0}&&{0}\\
\enddiagram
$$}
\noindent
With respect to suitably chosen bases of $U$ and $\Lambda^1 V^*\oplus
\Lambda ^0 V^*$ the map $\phi$ is given by the matrix
$$
\phi =
\pmatrix
 x_4&   0&   0&   0& x_1& x_3&   0&   0& x_2&   0&   0&   0&   0&   0\\
   0& x_4&   0&   0& -x_0&   0&x_2&   0&   0& x_3&   0&   0&   0&   0\\
   0&   0& x_4&   0&   0&   0&-x_1& x_3& -x_0&   0&   0&   0&   0&   0\\
   0&   0&   0& x_4&   0& -x_0&   0&-x_2&   0&-x_1&   0&   0&   0&   0\\
-x_0&-x_1&-x_2&-x_3&   0&   0&   0&   0&   0&   0&   0&   0&   0&   0\\
   0&   0&   0&   0&   0&   0&   0&   0& x_4& x_4& x_0& x_1& x_2& x_3  
\endpmatrix
\ 
.
$$
Resolving $\phi$ shows that $\scr K$ has syzygies of type
$$
\vbox{%
\halign{&\hfil\,$#$\,\hfil\cr
0&\leftarrow&\scr K&\leftarrow&&15 \scr O(-1)&&5\scr O(-2)&&\scr O(-3)\cr
&&&&&&\vbox to 10pt{\vskip-0pt\hbox{$\nwarrow$}\vss}&\oplus&\leftarrow&\oplus\cr
&&&&&&&5\scr O(-3) &&4\scr O(-4)&\vbox to 10pt{\vskip-4pt\hbox{$\nwarrow$}\vss}&
\scr O(-5)&\leftarrow&0\ .\cr
}}
$$
Comparing with the Koszul resolution of $\Omega^3(3)$ gives the result.
\qed
\enddemo
\proclaim
{Remark 2.3} 
Conversely, by starting with $(M)$ and $\scr G$, we obtain
a proof for the existence of smooth surfaces in $\openP^4$
with $d=8$, $\pi = 5$ and $q=1:$ Since $\scr G(1)$ is globally generated, 
Kleiman's Bertini-type theorem \cite {Kl} implies that the dependancy locus of 4 general 
sections of $\scr G(1)$ is indeed a smooth surface (with invariants as above). \quad\qed
\endproclaim
\proclaim
{Remark 2.4} (i) Let $X\overset \Phi\to\rightarrow C$ be an 
elliptic conic bundle as above. Ionescu \cite{Io} studied the 
union $W$ of the planes of the conics which are the fibres of $\Phi$. 
He showed, that $W$ is a quartic hypersurface in $\openP^4$ which 
is a cone with vertex a point. The base $S$ of the cone is an elliptic 
ruled quartic surface in $\openP^3$ with two double lines $l_1$
and $l_2$. The rulings of $S$ form a divisor of class (2,2) 
in the $\openP^1 \times \openP^1\cong l_1 \times l_2$ 
of lines joining $l_1$ and $l_2$. Ionescu considered the cone
over a pair of rulings through a point on $l_1$ but overlooked
the presence of $l_2$. His residual quadric at the end of 
\cite{Io, 6.4} is the double plane over $l_2$.\par\noindent
(ii) The vertex of the cone $W$ is the distinguished point 
$\langle e_4\rangle$ in our construction. The projection from
$\langle e_4\rangle$ maps $X$ 2:1 to $S$. The conics on $X$
are mapped 2:1 onto the rulings of $S$.
\endproclaim 
\proclaim
{Remark 2.5} In the minimal free resolution of $\scr G$ (and thus also of $\scr J_X$)
there are two maps which are given by four linear forms. In fact, for both maps, these 
forms generate the ideal of the distinguished point $\langle e_4\rangle$.\quad\qed
\endproclaim
\bigskip\noindent
\head
{3. Vector bundles via stratifications}
\endhead
In this section we take a quick look at the approach of the third author 
to the construction of vector bundles on $\openP^n$. In particular, we 
construct a stable rank-3 bundle $\scr E$ on $\openP^4$ such that
the dependancy locus of two general sections of $\scr E$ is an elliptic
conic bundle. This rank-3 bundle was the starting point of our
investigations. It is related to the rank-5 bundle $\scr G$ from
Section 2 via an exact sequence
$$
0\rightarrow 2 \scr O\rightarrow\scr G(1) \rightarrow \scr E\rightarrow 0.
$$
\medskip
Let $\scr E$ be a rank-$r$ vector bundle on $\openP^n$ with first Chern 
class $c_1$. Minimal systems of generators $s_1,\dots,s_l$ of $H^0_*\scr E$ 
and $\sigma_1,\dots, \sigma_k$ of $H^0_*\scr E^{\vee}$ give rise to a 
commutative diagram
$$
\diagram[height=1.5em]
\scr L&&\rTo^{S}&&{\scr K}\\
&\rdTo&&\ruTo\\
&&\scr E&\\
&0\ruTo(1,1)&&\quad \rdTo(1,1) 0\quad ,\\
\enddiagram
$$
where $\scr L$ and $\scr K$ are direct sums of line bundles of ranks $l$
and $k$ resp. $S$ is a $k\times l$-matrix with polynomial entries. The 
entries in the $j$-th column of $S$ define the zero-scheme $X_{s_j}=\{s_j=0\}$.
The entries in the $i$-th row of $S$ define the zero-scheme $X_{\sigma_i}
=\{\sigma_i=0\}$. Suppose that $\scr K$ has a direct summand of type 
$r\scr O(m)$, and let $\sigma_{i_1},\dots,\sigma_{i_r}$ be the corresponding 
elements of $H^0_*\scr E^{\vee}$. By projecting onto this direct summand we may
represent $\scr E$ as a subsheaf of $r\scr O(m)$:
$$
\diagram[height=1.5em]
\scr L&&\rTo^{T}&&{r\scr O(m)}\\
&\rdTo&&\ruTo_{\theta}\\
&&\scr E&\\
&0\ruTo(1,1)&&\quad \rdTo(1,1) 0\quad .\\
\enddiagram
$$
In fact, $\theta$ is a monomorphism of vector bundles outside the divisor defined 
by the form $f=\sigma_{i_1}\wedge\dots\wedge \sigma_{i_r}\in H^0(\openP^n,
\scr O((r\cdot m-c_1))\cong H^0(\openP^n,\Lambda^r \scr E^{\vee}(m))$. 
The $j$-th column of $T$ represents the section $t_j:=H^0(\theta)(s_j)$ 
of $r\scr O(m)$, and we have the relations of forms $t_{j_1}\wedge\dots
\wedge t_{j_r}=f\cdot (s_{j_1}\wedge\dots\wedge s_{j_r})$. In order
to detect properties of $\scr E$ it is often enough to study the matrix $T$.
For example, one may verify that $\scr E$ is indeed a vector bundle by looking
at the ideal $I$ generated by the maximal minors of $T$ and by checking
that the ideal quotient $(I:f)$ defines the empty subset of $\openP^n$.
\medskip
Conversely, one can construct vector bundles $\scr E$ on $\openP^n$ by starting 
with convenient matrices $T$. In the examples below $f$ is of the form 
$(f_{\red})^{(r-1)}$, where $f_{\red}=x_{0}\cdot\dots\cdot x_{c_1}$ is a product of 
coordinates, and $T$ is of the form $T=\pmatrix T^{\prime}&T^{\prime\prime}
\endpmatrix:\scr L \rightarrow r\scr O(c_1)$, where $T^{\prime}$ is the 
appropriate identy 
matrix multiplied by $f_{\red}$. In each example the minimal free resolution
of $\scr E$ is obtained by resolving $T$ (use e.g. Macaulay \cite {BS}).
\medskip\noindent
{\it {Example 3.1}}. The Nullcorrelation bundle $\scr E$ on $\openP^3$
(compare e.g. \cite {Ba}).\par\noindent
In this case $r=2$, $c_1=2$, and we may choose $f=f_{\red}=x_0x_1$ and 
$$
T=
\pmatrix
x_0x_1&0&x_0x_3&x_1x_2&x_2x_3\\
0&x_0x_1&x_0^2&x_1^2&x_0x_2+x_1x_3\\
\endpmatrix
.
$$ 
Note that the sections of $\scr E^{\vee}(2)\cong \scr E$ corresponding to the 
rows of $T$ vanish along the u\-ni\-on of two skew lines and a double line 
resp. The minimal free resolution of $\scr E$ is of type
$$
\vbox{%
\halign{&\hfil\,$#$\,\hfil\cr
0\leftarrow\scr E&\leftarrow5 \scr O\leftarrow 4\scr O(-1)\leftarrow
\scr O(-2)\leftarrow 0\cr
}}.\quad\qed
$$
\medskip\noindent
{\it {Example 3.2}}. The Horrocks-Mumford bundle $\scr E$ on $\openP^4$
\cite {HM}.
\par\noindent
In this case $r=2$, $c_1=5$, and we may choose $f=f_{\red}=x_0x_1x_2x_3x_4$
and 
$$
T=
\pmatrix
T_0&T_1&T_2&T_3\\
\endpmatrix,
$$ 
with $T_0,\dots ,T_3$ as follows:
$$
T_0=
\pmatrix
\gamma_0&0&\gamma_3&\gamma_4 \\
0&\gamma_0&\gamma_1&\gamma_2 \\
\endpmatrix ,
$$
where $\gamma_0,\dots ,\gamma_4$ are the five Horrocks-Mumford quintics
$$
\matrix
\gamma_0=y_0y_1y_2y_3y_4 &
\gamma_1=\sum_{i\in\openZ_5}y_iy_{i+2}^2y_{i+3}^2 &
\gamma_2=\sum_{i\in\openZ_5}y_i^3y_{i+2}y_{i+3} \\
\gamma_3=\sum_{i\in\openZ_5}y_i^3y_{i+1}y_{i+4} &
\gamma_4=\sum_{i\in\openZ_5}y_iy_{i+1}^2y_{i+4}^2 .
\endmatrix
$$
The entries of $T_1,\dots ,T_3$ are sextics:
$$
T_1=\pmatrix
x_{i+1}x_{i+2}x_{i+3}~^t(x_{i+2}x_{i+3}x_{i+4}\quad x_{i+1}x_{i+4}^2)
\endpmatrix _{i\in\openZ _5} ,
$$
$$
T_2=\pmatrix
x_{i+2}x_{i+3}x_{i+4}~^t(x_{i+2}x_{i+3}x_{i+4}\quad x_{i+1}x_{i+4}^2
+x_{i}x_{i+2}^2)\endpmatrix _{i\in\openZ _5},
$$
$$
T_3=\pmatrix
x_{i+1}x_{i+3}x_{i+4}~^t(x_{i+2}x_{i+3}^2+x_{i}x_{i+4}^2\quad 
x_{i+1}x_{i+3}x_{i+4})\endpmatrix _{i\in\openZ _5} .
$$
\par\noindent
The zero-schemes of the sections of $\scr E$ are called {\it {Horrocks-Mumford
surfaces}}. They are the minimal abelian surfaces in $\openP^4$ and the 
degenerations of these smooth surfaces \cite {HM} (compare \cite {BHM}). Those 
sections of $\scr E^{\vee}(5)\cong \scr E$ corresponding to the rows of $T$ vanish 
along the unions of five double planes, namely 
$\lbrace x_{i}^2=x_{i+2}^2=x_{i}x_{i+2}=x_{i+2}x_{i+3}^2+x_{i}x_{i+4}^2=0\rbrace$,
$i\in\openZ_5$, and 
$\lbrace x_{i}^2=x_{i+1}^2=x_{i}x_{i+1}=x_{i+1}x_{i+4}^2+x_{i}x_{i+2}^2=0\rbrace$,
$i\in\openZ_5$, resp. (compare e.g. \cite{ADHPR, Section 9} for these degenerations). 
The minimal free resolution of $\scr E$ is of type
$$
\vbox{%
\halign{&\hfil\,$#$\,\hfil\cr
&&&&4\scr O\cr
0&\leftarrow&\scr E&\leftarrow&\oplus\cr
&&&&15\scr O(-1)&\vbox to
10pt{\vskip-4pt\hbox{$\nwarrow$}\vss}&35\scr O(-2)&\leftarrow&20\scr
O(-3)\cr
&&&&&&&&&\vbox to
10pt{\vskip-4pt\hbox{$\nwarrow$}\vss}&2\scr O(-5)&\leftarrow &0\ \cr
}}
$$
(compare e.g. \cite{De}).\quad\qed
\medskip\noindent
{\it {Example 3.3}}. A stable rank-2 reflexive sheaf $\scr E$ on $\openP^4$
such that $\scr E(-3)$ has Chern classes $c_1=-1$, $c_2=9$, $c_3=25$ 
and $c_4=50$ \cite {ADHPR, Section 7}.
\par\noindent
In this case $r=2$, $c_1=5$, and we may choose $f=f_{\red}=x_0x_1x_2x_3x_4$
and 
$$
T=
\pmatrix
T_0&T_1\\
\endpmatrix,
$$ 
with $T_0$ and $T_1$ as follows:
$$
T_0=
\pmatrix
\gamma_0&0&\gamma_2&\gamma_4 \\
0&\gamma_0&\gamma_1&\gamma_3 \\
\endpmatrix ,
$$
where $\gamma_0,\dots ,\gamma_4$ are the Horrocks-Mumford quintics from 3.2.
The entries of $T_1$ are septics:
$$
T_1=\pmatrix
x_{i+1}x_{i+2}x_{i+3}x_{i+4}~^t(x_{i+2}x_{i+4}^2+x_{i+1}^2x_{i+3}\quad
x_{i+3}^2x_{i+4}+x_{i+1}x_{i+2}^2)\endpmatrix _{i\in\openZ _5}.
$$
\par\noindent
This time $\scr E$ is locally free only outside the union of the
25 Horrocks-Mumford lines (compare \cite {HM} for these lines). 
The zero-schemes of the sections of $\scr E$ are non-minimal abelian 
surfaces in $\openP^4$ and the degenerations of these smooth surfaces.
They are (5,5)-linked to Horrocks-Mumford surfaces. The zero-schemes of those 
sections of $\scr E^{\vee}(5)\cong \scr E$ corresponding to the rows of $T$ 
consist of five components, namely 
$\lbrace x_{i}=x_{i+1}^2x_{i+3}=x_{i+2}x_{i+4}^2=0\rbrace$, $i\in\openZ_5$, and 
$\lbrace x_{i}=x_{i+1}x_{i+2}^2+x_{i+3}^2x_{i+4}=0\rbrace$, $i\in\openZ_5$, resp.
The minimal free resolution of $\scr E$ is of type
$$
\vbox{%
\halign{&\hfil\,$#$\,\hfil\cr
&&&&4\scr O\cr
0&\leftarrow&\scr E&\leftarrow&\oplus\cr
&&&&5\scr O(-2)&\vbox to 10pt{\vskip-4pt\hbox{$\nwarrow$}\vss}&15\scr O(-3)&\leftarrow&10\scr
O(-4)&\leftarrow&2\scr O(-5)\leftarrow 0.\ \cr
}}
$$
Compare \cite{ADHPR, Section 7 and 9}.\quad\qed
\medskip\noindent
{\it {Example 3.4}}. A stable rank-3 vector bundle $\scr E$ on $\openP^4$
with Chern classes $c_1=4$, $c_2=8$ and $c_3=8$.
\par\noindent
We choose $f=f_{\red}=x_0x_1x_2x_3$
and 
$$
T=
\pmatrix
T^{\prime}&t_4&&t_5&&t_6&&t_7&&t_8\\
\endpmatrix,
$$ 
with
$$
T^{\prime}=
\pmatrix
x_0x_1x_2x_3&0&0 \\
0&x_0x_1x_2x_3&0 \\
0&0&x_0x_1x_2x_3 
\endpmatrix ,
$$
and
$$
t_4=\matrix
x_{0}x_{3}~^t(x_{4}^2+x_{1}x_{3}-x_{0}x_{2}\quad x_{2}x_{3}-x_{0}x_{1}\quad 
-x_{0}x_{3})
\endmatrix ,
$$
$$
t_5=\matrix
x_{0}x_{1}~^t(x_{4}^2+x_{1}x_{3}+x_{0}x_{2}\quad -x_{0}x_{1}\quad x_{1}x_{2}
-x_{0}x_{3})
\endmatrix ,
$$
$$
t_6=\matrix
x_{1}x_{2}~^t(-x_{4}^2+x_{1}x_{3}-x_{0}x_{2}\quad x_{0}x_{1}-x_{2}x_{3}\quad 
-x_{1}x_{2})
\endmatrix ,
$$
$$
t_7=\matrix
x_{2}x_{3}~^t(-x_{4}^2+x_{1}x_{3}+x_{0}x_{2}\quad -x_{2}x_{3}\quad x_{0}x_{3}
-x_{1}x_{2})
\endmatrix ,
$$
$$
t_8=\pmatrix
-x_{4}^4+x_{1}^2x_{3}^2+x_{0}^2x_{2}^2\\
x_{4}^2(x_{0}x_{1}-x_{2}x_{3})-(x_{0}x_{1}+x_{2}x_{3})(x_{1}x_{3}+x_{0}x_{2})\\ 
(x_{0}x_{2}-x_{1}x_{3})(x_{1}x_{2}+x_{0}x_{3})+x_{4}^2(x_{0}x_{3}-x_{1}x_{2})
\endpmatrix .
$$
\par\noindent
A check on the maximal minors of $T$ shows that $\scr E$ is locally free
(use e.g. Macaulay \cite {BS}). Let $X_1$, $X_2$ and $X_3$ be 
the zero-schemes of the sections of  $\scr E^{\vee}(4)$ 
corresponding to the rows of $T$. Then $X_1$ is the union of four plane quadrics, 
whereas $X_2$ and $X_3$ are the unions of four double planes (use e.g. the primary
decomposition package of Singular \cite{GPS} to obtain the explicit equations).
The minimal free resolution of $\scr E$ is of type
$$
\vbox{%
\halign{&\hfil\,$#$\,\hfil\cr
0&\leftarrow&\scr E&\leftarrow&&8 \scr O&&4\scr O(-1)&&\scr O(-2)\cr
&&&&&&\vbox to 10pt{\vskip-0pt\hbox{$\nwarrow$}\vss}&\oplus&\leftarrow&\oplus\cr
&&&&&&&5\scr O(-2) &&4\scr O(-3)&\vbox to 10pt{\vskip-4pt\hbox{$\nwarrow$}\vss}&
\scr O(-4)&\leftarrow&0\ .\cr
}}
$$
The minimal free resolution of $\scr E^{\vee}$ is of type
$$
\vbox{%
\halign{&\hfil\,$#$\,\hfil\cr
&&&&\scr O(-2)\cr
&&&&\oplus\cr
0&\leftarrow&\scr E^{\vee}&\leftarrow&4\scr O(-3)&&2\scr O(-4)\cr
&&&&\oplus&&\oplus\cr
&&&&8\scr O(-4)&\vbox to 10pt{\vskip-4pt\hbox{$\nwarrow$}\vss}&12\scr O(-5)&
&4\scr O(-6)\cr
&&&&&&\oplus&\vbox to 10pt{\vskip-4pt\hbox{$\nwarrow$}\vss}&\oplus\cr
&&&&&&5\scr O(-6)&&8\scr O(-7)&\vbox to 10pt{\vskip-4pt\hbox{$\nwarrow$}\vss}&
&3\scr O(-8)\leftarrow 0.\ \cr
}}
$$
In particular, $\scr E$ is stable since the normalized bundle $\scr E(-2)$ and its 
dual twisted by -1 have no non-zero sections.\par\noindent
Let $(M)$ be the monad of the rank-5 vector bundle $\scr G$ as in Proposition 2.2.
By taking sections in the display of $(M)$ we identify $H^0 \scr G(1)$ in terms of 
the exterior algebra $\Lambda V^*$ as the kernel
$$
0\rightarrow H^0 \scr G(1)\rightarrow \Lambda^1 V^*\oplus
\Lambda ^0 V^* \overset\beta\to\rightarrow V^* \rightarrow 0.
$$
Two general elements of $H^0 \scr G(1)$ are nowhere dependant. In fact, 
$$
\sigma_1=
\pmatrix
e_0\wedge e_1 + e_2\wedge e_3\\
0
\endpmatrix
$$
and
$$
\sigma_2=
\pmatrix
e_1\wedge e_2-e_0\wedge e_3\\
e_1\wedge e_2\wedge e_4 - e_0\wedge e_3\wedge e_4
\endpmatrix
$$
give rise to an exact sequence
$$
0\rightarrow 2 \scr O\overset\sigma\to\rightarrow\scr G(1) 
\rightarrow \scr E\rightarrow 0.
$$
This can be easily seen by comparing the syzygies of the cokernel of $\sigma$ with
those of $\scr E$ (use e.g. Macaulay \cite {BS}). It follows that $\scr E(-1)$ is 
isomorphic to the cohomology bundle of the monad
$$
\qquad 0\rightarrow 2\scr O(-1) \oplus \Omega ^3(3) 
\overset(\sigma,\alpha)\to\longrightarrow
\Omega ^2(2) \oplus \Omega ^1(1) \overset\beta\to\longrightarrow 
\scr O\rightarrow 0. \quad\qed
$$
In each of these examples we have a stratification
$$
L^{m_0} \subset\dots \subset L^1 \subset \openP^n,\quad m_0={\text{min}}
\{n, c_1\},
$$ 
where we use the following notations: For any subset $I$ of 
$\openZ_{m_0}$ we set
$L_I = \bigcap_{i\in I} \ L_i$, where $L_i=\lbrace x_i=0\rbrace$, 
and for any $m\in\openZ_{m_0}$ we write
$$
L^m=\bigcup \{L_I \mid I\in \Cal P_m(\openZ_{m_0})\}, 
$$
where $\Cal P_m(\openZ_{m_0})$ is the collection of subsets $I\subset 
\openZ_{m_0}$ whith cardinality $\mid I \mid =m$. One may use this
filtration to find local frames for $\scr E$ in a more systematic way,
and to construct a filtration of sheaves
$$
\scr E^0\cong r\scr O\subset\scr E^1\dots\subset\scr E^{m_0}\cong \scr E,
$$
such that
$$
\scr E^{m}/\scr E^{m-1}\cong \bigoplus_{I\in \Cal P_m(\openZ_{m_0})} 
\scr O_{L_I}(d_{L_I}), 
$$
with convenient twists $d_{L_I}$. Via this filtration one can e.g. compute the 
cohomology of $\scr E$ (use $\check{\C}$ech cohomology). We refer to 
\cite{Sa} and \cite{SEKS} for more details.
%

%
\enddocument